# On Berry–Esseen bounds for non-instantaneous filters of linear processes

TSUNG-LIN CHENG[1] and HWAI-CHUNG HO[2]

[1]*Department of Mathematics, National Changhua University of Education, Changhua, Taiwan. E-mail: tlcheng@math.ncue.edu.tw*

[2]*Institute of Statistical Science, Academia Sinica, and Department of Finance, National Taiwan University, Taipei 115, Taiwan. E-mail: hcho@stat.sinica.edu.tw*

Let $X_n = \sum_{i=1}^{\infty} a_i \varepsilon_{n-i}$, where the $\varepsilon_i$ are i.i.d. with mean 0 and at least finite second moment, and the $a_i$ are assumed to satisfy $|a_i| = O(i^{-\beta})$ with $\beta > 1/2$. When $1/2 < \beta < 1$, $X_n$ is usually called a *long-range dependent* or *long-memory* process. For a certain class of Borel functions $K(x_1, \ldots, x_{d+1})$, $d \geq 0$, from $\mathcal{R}^{d+1}$ to $\mathcal{R}$, which includes indicator functions and polynomials, the stationary sequence $K(X_n, X_{n+1}, \ldots, X_{n+d})$ is considered. By developing a finite orthogonal expansion of $K(X_n, \ldots, X_{n+d})$, the Berry–Esseen type bounds for the normalized sum $Q_N/\sqrt{N}, Q_N = \sum_{n=1}^{N}(K(X_n, \ldots, X_{n+d}) - EK(X_n, \ldots, X_{n+d}))$ are obtained when $Q_N/\sqrt{N}$ obeys the central limit theorem with positive limiting variance.

## 1. Introduction

Consider a linear process $X_n$, $n = 1, 2, \ldots$, defined by $X_n = \sum_{i=1}^{\infty} a_i \varepsilon_{n-i}$, where the $\varepsilon_i$ are i.i.d. having mean 0 and at least finite second moment and the $a_i$ are assumed to satisfy $|a_i| = O(i^{-\beta})$ with $\beta > 1/2$. Let $K(x_1, \ldots, x_{d+1})$, $d \geq 0$, be a Borel function from $\mathcal{R}^{d+1}$ to $\mathcal{R}$. For fixed $d \geq 0$, define $Q_N = \sum_{n=1}^{N}(K(X_n, \ldots, X_{n+d}) - EK(X_n, \ldots, X_{n+d}))$. The present paper aims to establish the Berry–Esseen-type of rate of convergence for $Q_N/\sqrt{N}$ when $Q_N/\sqrt{N}$ obeys the central limit theorem with some positive limiting variance $\sigma^2$, that is, to determine a positive real number $s$ such that the uniform distance $\sup_x |P(Q_N/\sqrt{N} \leq x) - \Phi(x/\sigma)|$ between the two distributions $P(Q_N/\sqrt{N} \leq x)$ and $\Phi(x/\sigma)$ is $O(N^{-s})$ as $N$ tends to infinity, where $\Phi(\cdot)$ is the standard normal distribution function. There is abundant literature investigating the same problem for i.i.d. sequences or stationary sequences which are weakly dependent (short-range dependent or short-memory; see Bradley [4] for a review on sequences of weak dependence). For the former, comprehensive studies are given in, for example, Gnedenko and Kolmogorov [7] and Petrov [18]. For the latter, which includes the case of $m$-dependence, a detailed account of related results can be found in Sunklodas [22]. With an emphasis on the statistical inference front, Bentkus, Götze and Tikhomirov [1] studied the Berry–Esseen bounds for a general class of asymptotically normal statistics constructed from absolutely regular random variables. In addition, considering a Gaussian linear process, Taniguchi [24] studied a problem similar to ours and derived a bound of order $O(N^{-1/2})$ by assuming that $\beta > 2$ and that the functional $K$ is of the form $K(x, y) = xy$. There are three features that distinguish our setting from those in the literature mentioned above. First, the function $K$ need not be smooth as considered in Taniguchi [24]. Second, the sequence $\{X_n\}, n \geq 1$, is not necessarily Gaussian. Third, the dependence structure of $\{X_n\}$ is determined solely through the decay rate $|a_i| = O(i^{-\beta})$ of the innovation coefficients without assuming any mixing-type condition, which is in general, difficult to verify. Moreover, when $\beta \in (1/2, 1)$, $\{X_n\}$ covers a widely studied class of long-range dependent (or long-memory) processes (cf. Brockwell and Davis [3]) which is not strong mixing (Rosenblatt [21]). Note that the mixing-type conditions are not satisfied by long-memory processes with $\beta < 1$, but also sometimes fail to hold even when the hyperbolic decay rate $\beta$ is greater than 1 (see Bradley [4], Phan and Tran [19]). With the same setting of non-instantaneous functionals and linear processes, but mainly dealing with the short-memory case were $\beta > 1$, Wu [27] proved central limit theorems for $Q_N/\sqrt{N}$ by treating the innovation sequence $\{\ldots, \varepsilon_{n-1}, \varepsilon_n\}$ as an infinite-dimensional Markov chain. When $\{X_n\}$ is long-range dependent with $1/2 < \beta < 1$, it is important to note that the asymptotic behavior of $Q_N$ depends very much on the function $K$. Although the variance $\text{var}(\sum_{n=1}^{N} X_n)$ of partial sums of $X_n$'s grows with the rate $O(N^{-2\beta+3})$, it is possible that $\text{var}(Q_N) = O(N)$ for certain $K$'s and that the root-$N$ central limit theorem holds. To prove this central limit theorem for $Q_N/\sqrt{N}$, the standard approach is to expand $K(X_n, \ldots, X_{n+d+1})$ in terms of polynomials, provided that either $\{X_n\}$ is Gaussian (Ho and Sun [15]) or the functional $K$ under consideration is smooth (Giraitis [8]). To deal with the case where both assumptions fail, Ho and Hsing [14] introduced a new method based on martingale decomposition to prove the central limit theorem for instantaneous functionals which include indicator functions and polynomials. In order to refine the central limit theorem by giving Berry–Esseen bounds for the rate of convergence, we shall combine the martingale method proposed in Ho and Hsing [14] with the blocking method (Bernstein [2]). The latter has been extensively used in studying the asymptotic behavior of the sum of weakly dependent random variables. It is known that under the weak dependence, the rates achieved by using the blocking method is not as sharp as those achieved by the Stein–Tikhomirov method (Stein [23], Tikhomirov [26]) which involves using a linear differential equation in terms of the difference between the distribution (characteristic) function of the sum of weakly dependent random variables and that of a standard normal random variable. However the blocking method is still an appealing technique in our setting since the Stein–Tikhomirov approach is hard to apply in the case where the condition of weak dependence no longer holds.

The main contribution of the present paper is the derivation of the Berry–Esseen bounds (given in (5) and (6) in Section 2.2) for the class of stationary sequences generated by applying a nonlinear transformation to linear processes that are allowed to be long-range dependent (when $1/2 < \beta < 1$) or short-range dependent (when $1 < \beta$), but need not satisfy any traditional mixing-type condition. The bounds will depend on the functional $K$ as well as on the decay rate $\beta$ of the coefficients $\{a_i\}$. It is known that the Berry–Esseen rate established for independent or $m$-dependent sequences is of order $O(N^{-1/2})$ (Gnedenko and Kolmogorov [7], Petrov [18] and Stein [23]). For weakly dependent sequences such as sequences that are strong mixing, absolutely regular or completely regular, the rate is bounded by $O(N^{-1/2} \log^2 N)$ provided that the mixing coefficients decay exponentially (Tikhomirov [26] and Stein [23]). While previous results in the literature require that the sequences be independent or weakly dependent, we demonstrate for the stationary sequence $\{K(X_n, \ldots, X_{n+d})\}$ that the Berry–Esseen bounds can still be obtained, even when the underlying process is long-range dependent. Although the rate we achieve is slower than $O(N^{-1/2})$, it is as good as those shown for a certain class of strong mixing sequences whose mixing coefficients decrease hyperbolically to zero (Theorem 1, Tikhomirov [26]).

The rest of this paper is organized as follows. In the next section, following the introduction of some notation and technical conditions, are the statements of our major results, Theorems 1 and 2, which deal with the long- and short-memory cases, respectively. Toward the end of Section 2, we present two commonly seen examples, zero crossings and lag covariances, to illustrate the results. Section 3 presents the proofs of the theorems, which make frequent use of two technical lemmas that are presented in Section 4.

## 2. Main results

Before stating our major results (Theorems 1 and 2 and Corollary 1 in Section 2.2), we give a brief description of our approach. Define $Z_n = K(X_n, X_{n+1}, \ldots, X_{n+d}) - \mathrm{E}K(X_n, X_{n+1}, \ldots, X_{n+d})$ and $X_{n,m} = \sum_{i=1}^{m} a_i \varepsilon_{n-i}$. We first extend the martingale decomposition technique introduced in Ho and Hsing [14] to the case of non-instantaneous functionals to establish that there exist an integer $r \geq 1$ and a real sequence $\{b_{j_1, \ldots, j_r}\}$ such that the sequence of autocovariances of the sequence $\{Z_n\}$ decays with the same rate as that of $\sum_{d \leq j_1, \ldots, j_r < \infty} b_{j_1, \ldots, j_r} \prod_{s=1}^{t+1} \varepsilon_{n+d-j_s}$. Here, the integer $r$ is determined by the functional $K$ and the underlying sequence $\{X_n\}$, and the coefficient $b_{j_1, \ldots, j_r}$ is asymptotically of order $O(j_1^{-\beta} \cdots j_r^{-\beta})$. Suppose that $1/2 < \beta < 1$ (the long-memory case) and $r(2\beta - 1) > 1$. It is then clear that the sequence of autocovariances of $\{Z_n\}$ is summable, that is, $\sum_{k=0}^{\infty} |\mathrm{E}Z_n Z_{n+k}| < \infty$. This property suggests that even when $\{X_n\}$ is long-memory, one can apply a certain nonlinear transformation to it to obtain a sequence $\{Z_n\}$ which behaves like a short-memory one and which obeys the central limit theorem. Based on this, it is then plausible to approximate $Z_n$ by $Z_{n,\ell(N))}$, with $Z_{n,\ell(N)} = K(X_{n,\ell(N)}, X_{n+1,\ell(N)}, \ldots, X_{n+d,\ell(N)}) - \mathrm{E}K(X_{n,\ell(N)}, X_{n+1,\ell(N)}, \ldots, X_{n+d,\ell(N)})$, where $\ell(N)$ increases to infinity at an appropriate rate as $N \to \infty$. For fixed $N$ and suitably chosen increasing sequences $\{\ell(N)\}, \{k_N\}$ and $\{A_N\}$, we adopt the blocking method (Bernstein [2]) to select $k_N$ disjoint blocks of $Z_{n,\ell(N)}$'s from $\{Z_{1,\ell(N)}, \ldots, Z_{N,\ell(N)}\}$, each having size $A_N$, such that these blocks are sufficiently far apart from each other and thus mutually independent. Our results then follow from deriving the Berry–Esseen bound for the $k_N$ independent blocks and letting the $N$ tend to infinity. The specific form of $b_{j_1,\ldots,j_r}$'s and the precise values of $k_N$, $\ell_N$ and $A_N$ will be given in the next subsection and during the course of the proofs, respectively.

### 2.1. Notation and technical conditions

For $u \geq 1$ define $(d+1)$-dimensional vectors of partial differentiation operators as

$$\mathbf{A}_u = \left(\mathbf{A}_{u,1}\frac{\partial}{\partial x_1}, \ldots, \mathbf{A}_{u,1+d}\frac{\partial}{\partial x_{d+1}}\right),$$

where

$$\mathbf{A}_{u,v} = \begin{cases} a_{u-(d+1-v)}, & \text{if } u-(d+1-v) > 0, \\ 0, & \text{if } u-(d+1-v) \leq 0, \end{cases}$$

for $1 \leq v \leq d+1$. For example,

$$\mathbf{A}_1 = \left(0, \ldots, 0, a_1\frac{\partial}{\partial x_{d+1}}\right),$$

$$\mathbf{A}_2 = \left(0, \ldots, 0, a_1\frac{\partial}{\partial x_d}, \ldots, a_2\frac{\partial}{\partial x_{d+1}}\right), \ldots$$

and for $u \geq d+1$,

$$\mathbf{A}_u = \left(a_{u-d}\frac{\partial}{\partial x_1}, \ldots, a_u\frac{\partial}{\partial x_{d+1}}\right).$$

Also, define

$$\mathbf{B}_j = \sum_{i=1}^{d+1} \mathbf{A}_{j,i}\frac{\partial}{\partial x_i}$$

and for $r \geq 2$,

$$\mathbf{B}_{j_1\cdots j_r} = \sum_{u_1,\ldots,u_r=1}^{d+1} \mathbf{A}_{j_1,u_1}\cdots\mathbf{A}_{j_r,u_r}\frac{\partial^r}{\partial x_{u_1}\cdots\partial x_{u_r}},$$

where $\partial^r/\partial x_{u_1}\cdots\partial x_{u_r}$ denotes partial differentiation with respect to the variables $x_{u_1}, \ldots, x_{u_r}$, $r$ times. That is, for any smooth function $G(\cdot)$,

$$B_{j_1\cdots j_r} \circ G(x_1, \ldots, x_{d+1}) = \sum_{u_1,\ldots,u_r=1}^{d+1} \mathbf{A}_{j_1,u_1}\cdots\mathbf{A}_{j_r,u_r}\frac{\partial^r G(x_1,\ldots,x_{d+1})}{\partial x_{u_1}\cdots\partial x_{u_r}}. \tag{1}$$

Recall that $X_{n,j} = \sum_{i=1}^{j} a_i\varepsilon_{n-i}$. For $1 \leq j \leq \infty$, define $\tilde{X}_{n,j} = X_n - X_{n,j}$ and

$$\mathbf{X}_{n,j} = (X_{n,(j-d)\vee 0}, X_{n+1,(j-d+1)\vee 0}, \ldots, X_{n+d,j}),$$

$$\tilde{\mathbf{X}}_{n,j} = \mathbf{X}_n - \mathbf{X}_{n,j} = (\tilde{X}_{n,(j-d)\vee 0}, \tilde{X}_{n+1,(j-d+1)\vee 0}, \ldots, \tilde{X}_{n+d,j}).$$

Let

$$\mathbf{X}_n = (X_n, X_{n+1}, \ldots, X_{n+d}),$$
$$\mathbf{x} = (x_1, \ldots, x_{d+1}).$$

For $1 \leq j \leq \infty$ and fixed $d$, let $[j]$ denote the set of $d+1$ indices $(j-d) \vee 0, (j-d+1) \vee 0, \ldots, j-1, j$ and $F_{[j]}$ the joint distribution function of $\mathbf{X}_{1,j}$. Note that $\mathbf{X}_{n,\infty} = \mathbf{X}_n$, that $F_{[\infty]}$ denotes the joint distribution of $(X_1, \ldots, X_{d+1})$ and that if any index in $[j]$ is zero, its corresponding distribution is set to be a point mass at zero with probability one. Define

$$K_{[j]}(x_1, \ldots, x_{d+1}) = \int K(x_1 + u_1, \ldots, x_{d+1} + u_{d+1}) \, dF_{[j]}(u_1, \ldots, u_{d+1})$$

or, in abbreviated from,

$$K_{[j]}(\mathbf{x}) = \int K(\mathbf{x} + \mathbf{u}) \, dF_{[j]}(\mathbf{u}).$$

Let $\mathbf{i} = (i_1, \ldots, i_{d+1})$. We denote the $\mathbf{i}$th partial derivative of $K_{[j]}(x_1, \ldots, x_{d+1})$ by

$$K_{[j]}^{\mathbf{i}}(x_1, \ldots, x_{d+1}) = \frac{\partial^{i_1, \ldots, i_{d+1}} K_{[j]}(x_1, \ldots, x_{d+1})}{\partial x_1^{i_1} \cdots \partial x_{d+1}^{i_{d+1}}}$$

whenever it exists. Define $Z_{N,0} = N E K(\mathbf{X}_n)$ and for $r \geq 1$,

$$Z_{N,r} = \sum_{n=1}^{N} \sum_{1 \leq j_1 < \cdots < j_r < \infty} \mathbf{B}_{j_1 \cdots j_r} \circ K_\infty(\mathbf{0}) \prod_{s=1}^{r} \varepsilon_{n+d-j_s}$$

and

$$Q_{N,p} = \sum_{n=1}^{N} K(\mathbf{X}_n) - \sum_{r=0}^{p} Z_{N,r}, \qquad p \geq 0.$$

Note that for the instantaneous case, that is, $d = 0$, $\mathbf{B}_{j_1 \cdots j_r} \circ K_\infty(\mathbf{0}) = a_{j_1} \cdots a_{j_r} K_\infty(0)$ and $Z_{N,r} = K_\infty(0) \sum_{n=1}^{N} \sum_{1 \leq j_1 < \cdots < j_r < \infty} \prod_{s=1}^{r} a_{j_s} \varepsilon_{n-j_s}$, which is precisely $K_\infty(0) Y_{N,r}$ as defined in Ho and Hsing ([14], page 1638).

Below are some regularity conditions that will be needed for the results stated in the next subsection. Let $\mathbf{i} = (i_1, \ldots, i_{d+1})$ and $\mathbf{x} = (x_1, \ldots, x_{d+1})$.

(C1) For a certain positive integer $J$, the partial derivatives $K_{[d+1]}^{\mathbf{i}}(\mathbf{x})$ of $K_{[d+1]}(\mathbf{x})$ of order $\mathbf{i} = (i_1, \ldots, i_{d+1})$ with $0 \leq i_1 + \cdots + i_{d+1} \leq J+2$ are continuous and one of the following two conditions holds:

(i) $K_{[d+1]}^{\mathbf{i}}(\mathbf{x})$ is bounded and $E\varepsilon_1^8 < \infty$;

(ii) $K^{\mathbf{i}}_{[d+1]}(\mathbf{x})$ is unbounded, but there is a polynomial function $U_{\mathbf{i}}(\mathbf{x})$ of degree $M$ such that $|K^{\mathbf{i}}_{[d+1]}(\mathbf{x})| \leq |U_{\mathbf{i}}(\mathbf{x})|$ for all $\mathbf{x} \in \mathcal{R}^{d+1}$, and $\mathrm{E}\varepsilon_1^{\max\{8,4M\}} + \mathrm{E}U_{\mathbf{i}}^4(X_1) < \infty$.

(C2) $\mathrm{E}[K(\mathbf{X}_1) - K(\mathbf{X}_{1,\ell})]^2 \to 0$ as $\ell \to \infty$.

Condition (C1) that describes a concrete class of transformations $K$ is not presented in the full generality as given in Ho and Hsing [14], yet it covers most of interesting cases in the literature. We choose to use (C1) merely for presentational simplicity since our main purpose is to introduce a method to obtain a rate of convergence in the current setting rather than to seek a class of transformations $K$ as general as possible. Note that in part (i) of (C1), indicator functions are included if the distribution function $G$ of $\varepsilon_1$ is sufficiently smooth. Condition (C1) ensures the following useful property needed later for proving theorems: for $(i_1, \ldots, i_{d+1})$ with $0 \leq i_1 + \cdots + i_{d+1} \leq J+2$ and $j \geq d+1$, $K^{(i_1,\ldots,i_{d+1})}_{[j]}$ is continuous and satisfies

$$K^{(i_1,\ldots,i_{d+1})}_{[j]}(x_1,\ldots,x_{d+1}) = \int K^{(i_1,\ldots,i_{d+1})}_{[j-1]}(x_1+y_1,\ldots,x_{d+1}+y_{d+1})\,\mathrm{d}G_{j-d}(y_1)\cdots \mathrm{d}G_j(y_{d+1}), \tag{2}$$

where $G_u$ is the distribution of $a_u\varepsilon_1$. (2) can be shown by using an argument similar to that used in proving Lemma 2.1 of Ho and Hsing [14]. Condition (C2) is a technical assumption and seems to be a natural assumption for the $\ell$-truncation argument we employ in Section 3 for proving theorems.

### 2.2. Theorems

In characterizing the limiting theorems for the case of instantaneous transformations, Ho and Hsing [14] proposed a quantity called *power rank*, which is analogous to the Hermite rank when $X_t$ is Gaussian. In the following, the multivariate version of the power rank is introduced for non-instantaneous transformations.

*Definition 1.* We say that the $d+1$-dimensional transformation $K$ has power rank $\nu$ if all of the partial derivatives $K^{(i_1,\ldots,i_{d+1})}_{[\infty]}(0,\ldots,0)$ of $K_{[\infty]}(0,\ldots,0)$ of order $i_1,\ldots,i_{d+1}$ with $i_1 + \cdots + i_{d+1} = s \leq \nu$ exist and the following conditions are satisfied: $K^{(i_1,\ldots,i_{d+1})}_{[\infty]}(0,\ldots,0) = 0$ if $i_1 + \cdots + i_{d+1} \leq s < \nu$, and there exists $(i'_1,\ldots,i'_{d+1})$ with $i'_1 + \cdots + i'_{d+1} = \nu$ such that $K^{(i'_1,\ldots,i'_{d+1})}_{[\infty]}(0,\ldots,0)$ is non-zero.

**Theorem 1.** *Assume that $|a_i| = O(i^{-\beta})$, $\beta \in (1/2, 1)$, and that conditions (C1) and (C2) hold. Let $p$ be any positive integer satisfying $J \geq p+1 > (2\beta-1)^{-1}$. Then*

$$N^{-1/2}Q_{N,p} \xrightarrow{d} N(0,\sigma^2), \tag{3}$$

where $\sigma^2 = \lim_{N\to\infty} N^{-1}\mathrm{var}(Q_{N,p})$. *Assume, furthermore, that $\sigma^2 > 0$ and that for some $\alpha_1 > 0$*

$$\mathrm{E}[K(\mathbf{X}_1) - K(\mathbf{X}_{1,\ell})]^2 = O(\ell^{-\alpha_1}) \qquad as\ \ell \to \infty. \tag{4}$$

*Then*

$$\sup_x |P(N^{-1/2}Q_{N,p} \leq x) - \Phi(x/\sigma)| = O(N^{-Q'/(3(2Q'+1))}), \tag{5}$$

*where $Q' = \alpha_1 \wedge (2\beta - 1) \wedge ((p+1)(2\beta - 1) - 1)$.*

**Remark 1.** (i) For the case of long-range dependence, if the power rank is one (e.g., $K$ is linear and $p = 0$), then the formula for of $Q'$ given in Theorem 1 is not applicable. even though the limit is Gaussian. This is consistent with the fact that in the circumstances, the normalized partial sums follow a non-central limit theorem, that is, the normalization constant is $N^H$, with $H > 1/2$. For the general long-memory case where $(p+1)(2\beta - 1) < 1$ with $p \geq 1$, then one has non-Gaussian limiting distributions. Consider, for example, the case where $K$ is univariate and has rank equal to 2. The limit is then usually referred to as the *Rosenblatt distribution* (Taqqu [25]), having a close form of characteristic function. A Berry–Esseen bound in these circumstances is feasible to obtain and will be discussed in a future paper. For the general case of nonlinear and non-instantaneous $K$ and non-Gaussian $\{X_n\}$, however, the technical difficulties caused by lacking analytical forms of the characteristic functions of both the partial sums $N^{-H} \sum_{n=1}^N (K(X_n, \ldots, X_{n+d}) - \mu)$ and their limit still remain to be overcome. When $K(x) = x$ is the identity function and $X_n$ is linear without normality assumption, Hall [10] obtained the Edgeworth expansion for the sample mean. The result obtained relies on the fact that the limit is Gaussian and the characteristic function of the sample mean can be written down explicitly by utilizing the linear structure of $X_n$. (ii) The rate given in (5) can be slightly improved if condition (4) is strengthened to $\mathrm{E}[K(X_1, X_{1,\ell})]^{2k} = O(\ell^{-\alpha(k)})$ for some positive integer $k > 1$ (as suggested by the referee). In order to achieve the improvement, one needs to deal with the complicated expansion of $(Q_{N,p} - Q_{N,p,\ell})^{2k}$ (see (7) for the definition of $Q_{N,p,\ell}$) to compute its expectation. The details involved are tedious and highly technical, and are omitted in the present paper.

An important implication of Theorem 1, as stated in the following corollary, is the case where the functionals are of power rank greater than one. We first note that if the power rank of $K$ is $p + 1 \geq 2$, then for each $1 \leq r \leq p$, all $\mathbf{B}_{j_1 \cdots j_r} \circ K_\infty(\mathbf{0})$'s vanish and, as a result, $Z_{N,r} = 0$ and $Q_N = Q_{N,p}$.

**Corollary 1.** *Let $\beta \in (1/2, 1)$ and $p + 1$ be the power rank of $K$. Assume conditions (C1) and (C2) hold. If $J \geq p + 1 > (2\beta - 1)^{-1}$, then*

$$N^{-1/2} Q_N \xrightarrow{d} N(0, \sigma^2) \qquad as\ N \to \infty.$$

*Furthermore, the convergence rate (5) holds if (4) is satisfied.*

**Theorem 2.** *Assume* $|a_i| = O(i^{-\beta}), \beta > 1,$ *and that conditions (*C1*) and (*C2*) hold. Then*

$$N^{-1/2} Q_N \xrightarrow{d} N(0, \sigma^2) \qquad as\ N \to \infty. \tag{6}$$

*Under assumption (4), (5) holds for* $Q' = \alpha_1 \wedge (2\beta - 2)$.

**Remark 2.** For condition (4), if $K$ is a polynomial function of degree $D$, then $\alpha_1 = 2\beta - 1$ provided that $E\varepsilon_1^{2D} < \infty$.

**Remark 3.** (i) Note that in Theorem 2 that a Berry–Esseen bound can be achieved regardless of the power rank of the functionals. (ii) For the short-range dependent $\{X_n\}$ considered in Theorem 2, if $K(x) = x$, an $O(N^{-1/2})$ bound can be obtained by applying Theorem 2 stated in Dedecker and Prieur [6] under an additional dependence condition ((7.19) in Dedecker and Prieur [6]) which is in terms of some new mixing coefficients and weaker than the traditional mixing-type coefficients. The condition, however, would require stronger summability restrictions on the innovation coefficients $\{a_i\}$. It is still unknown whether the same $O(N^{-1/2})$ rate can be extended to the case of general functionals $K$.

**Remark 4.** We use two examples to further illustrate the theorems. Recall that our underlying process is defined as $X_n = \sum_{i=1}^{\infty} a_i \varepsilon_{n-i}$, where the $\varepsilon_i$ are mean zero i.i.d. and have at least finite second moment and where $|a_i| = O(i^{-\beta})$ with $\beta > 1/2$. (i) Zero crossings. Assume that $\{X_n\}$ is long-range dependent with $1/2 < \beta < 1$, and that the distribution of $X_n$ is standard normal. Define the functional $K$ as $K(X_n, X_{n+1}) = 1$ if $X_n X_{n+1} < 0$ and $0$ if $X_n X_{n+1} \geq 0$. In other words, $\sum_{n=1}^{N} K(X_n, X_{n+1})$ counts the total number of times that the sample path of $\{X_1, \ldots, X_N\}$ crosses zero. (See Kedem [16] for useful applications of zero crossings in time series analysis.) It is indicated in Ho and Hsing ([14], Remark 2 on page 1640) that when $\{X_n\}$ is Gaussian and the functional is instantaneous (i.e., $d = 0$), the power rank is identical to the Hermite rank as defined in Taqqu [25]. This property can be seen to hold for the non-instantaneous case. Therefore, the power rank of $K$ is two since its Hermite rank is, as computed in Ho and Sun [15], two. Suppose, further, that $2(2\beta - 1) > 1$. Then, by Corollary 1, the zero-crossing counts obey the central limit theorem with convergence rate as specified in (5), where $\alpha_1$ can be shown to be $2\beta - 1$ since $K(X_n, X_{n+1})$ is of power rank two. (ii) Lag covariances. The functional we consider is $K(X_n, \ldots, X_{n+d}) = X_n X_{n+d}$, which is frequently used when estimating lag covariances. While there has been much work done on central limit theorems on lag covariances for stationary sequences (see, e.g., Hannan [12], Hall and Heyde [11], Giraitis and Surgailis [9], Phillips and Solo [20], Hosking [13], Wu and Min [28], among others), the issue of Berry–Esseen bounds in the framework considered in Theorems 1 or 2 has not been addressed before. Straightforward computation shows that the power rank of $K$ is two. Corollary 1 can then be applied by the same argument as in the preceding example (1). Also, note that one can compute the value of $b_{j_1, j_2} \equiv \mathbf{B}_{j_1, j_2} \circ K_{[\infty]}(0, 0)$ by using formula (1), which can, in this particular case, be alternatively verified by simply multiplying $X_n X_{n+d}$ out and comparing the coefficients. For example, since $K(x_1, \ldots, x_{d+1}) = x_1 x_{d+1}$, it follows from formula (1) that $b_{1,j} = 0$ for $1 < j \leq d$ and $a_1^2$ for $j = d+1$, which coincide with the corresponding coefficients in the expansion of $X_n X_{n+d}$ by multiplication. More specifically, in the expansion of

$$X_n X_{n+d} - \mathrm{E} X_n X_{n+d} = \sum_{j \neq i+d} a_i a_j \varepsilon_{n-i} \varepsilon_{n+d-j} + \sum_{i=1}^{\infty} a_i a_{i+d}(\varepsilon_{n-i}^2 - \mathrm{E}\varepsilon_1^2),$$

the second (or square) term follows the $\sqrt{N}$ central limit theorem since $\sum_{i=1}^{\infty} |a_i a_{i+d}| < \infty$, while it requires the extra condition $2(2\beta - 1) > 1$ for the same asymptotic normality to hold for the first (or cross product) term.

## 3. Proofs

Recall that

$$Z_{N,r} = \sum_{n=1}^{N} \sum_{1 \leq j_1 < \cdots < j_r < \infty} \mathbf{B}_{j_1 \cdots j_r} \circ K_\infty(\mathbf{0}) \prod_{s=1}^{r} \varepsilon_{n+d-j_s},$$

$$Q_{N,p} = \sum_{n=1}^{N} K(\mathbf{X}_n) - \sum_{r=0}^{p} Z_{N,r}, \qquad p \geq 0.$$

Define the truncated versions of $Z_{N,r}$ and $Q_{N,p}$ as

$$Z_{N,r,\ell} = \sum_{n=1}^{N} \sum_{1 \leq i_1 < \cdots < i_r \leq \ell} \mathbf{B}_{j_1 \cdots j_r} \prod_{s=1}^{r} \varepsilon_{n+d-j_s}, \qquad r \geq 1,$$

and

$$\begin{aligned}
Q_{N,p,\ell} &= \sum_{n=1}^{N} K(\mathbf{X}_{n,\ell}) - \sum_{r=0}^{p} Z_{N,r,\ell}. \\
&\equiv \sum_{n=1}^{N} T_n(p,n).
\end{aligned} \qquad (7)$$

Also, define

$$\tilde{\mathbf{X}}_{n,j,\ell} = \tilde{\mathbf{X}}_{n,j} - \tilde{\mathbf{X}}_{n,\ell}, \qquad 1 \leq j \leq \ell.$$

The main building block of our proof is the following martingale decomposition of $K(\mathbf{X}_n) - \mathrm{E} K(\mathbf{X}_n)$:

$$K(\mathbf{X}_n) - \mathrm{E} K(\mathbf{X}_n) = \sum_{j=1}^{\infty} [K_{[j-1]}(\tilde{\mathbf{X}}_{n,j-1}) - K_{[j]}(\tilde{\mathbf{X}}_{n,j})], \qquad (8)$$

where $K_{[0]} = K$. Fix $n, n', j, j'$ and write

$$K_{[i]}(\tilde{\mathbf{X}}_{m,i}) = \mathrm{E}(K(\mathbf{X}_m)|\mathcal{F}_{m+d-i-1}),$$

where $\mathcal{F}_s$ is the $\sigma$-field generated by $\varepsilon_k, k \leq s$. Suppose that $n - j \neq n' - j'$ and, without loss of generality, assume that $n - j < n' - j'$. Then

$$\begin{aligned}
\mathrm{E}&[K_{[j-1]}(\tilde{\mathbf{X}}_{n,j-1}) - K_{[j]}(\tilde{\mathbf{X}}_{n,j})][K_{[j'-1]}(\tilde{\mathbf{X}}_{n',j'-1}) - K_{[j']}(\tilde{\mathbf{X}}_{n',j'})] \\
&= \mathrm{E}\{[\mathrm{E}(K(\mathbf{X}_n)|\mathcal{F}_{n+d-j}) - \mathrm{E}(K(\mathbf{X}_n)|\mathcal{F}_{n+d-j-1})] \\
&\quad \times [\mathrm{E}(K(\mathbf{X}_{n'})|\mathcal{F}_{n'+d-j'}) - \mathrm{E}(K(\mathbf{X}_{n'})|\mathcal{F}_{n'+d-j'-1})]\} \\
&= \mathrm{E}\{[\mathrm{E}(K(\mathbf{X}_n)|\mathcal{F}_{n+d-j}) - \mathrm{E}(K(\mathbf{X}_n)|\mathcal{F}_{n+d-j-1})] \\
&\quad \times \mathrm{E}[\mathrm{E}(K(\mathbf{X}_{n'})|\mathcal{F}_{n'+d-j'}) - \mathrm{E}(K(\mathbf{X}_{n'})|\mathcal{F}_{n'+d-j'-1})|\mathcal{F}_{n+d-j}]\} \\
&= \mathrm{E}\{[\mathrm{E}(K(\mathbf{X}_n)|\mathcal{F}_{n+d-j}) - \mathrm{E}(K(\mathbf{X}_n)|\mathcal{F}_{n+d-j-1})] \\
&\quad \times [\mathrm{E}(K(\mathbf{X}_{n'})|\mathcal{F}_{n+d-j}) - \mathrm{E}(K(\mathbf{X}_{n'})|\mathcal{F}_{n+d-j})]\} \\
&= 0.
\end{aligned} \qquad (9)$$

We now build a representation for $Q_{N,p} - Q_{N,p,\ell}$, which will be central to the proofs, based on the martingale decomposition (8). The main step to achieve the representation is to use $\sum_{r=1}^{p}(\prod_{s=1}^{r}\varepsilon_{n+d-j_s})\mathbf{B}_{j_1\cdots j_r} \circ K_{[\infty]}(\mathbf{0})$ for suitable $p$ to approximate the summand $K_{[j-1]}(\tilde{\mathbf{X}}_{n,j-1}) - K_{[j]}(\tilde{\mathbf{X}}_{n,j})$ (for $j \geq d+1$) by repeated applications of the martingale decomposition technique and differentiation. The task is carried out in a similar fashion for both the $Q_{N,p}$ and its truncated version $Q_{N,p,\ell}$. Write

$$Q_{N,p} - Q_{N,p,\ell} = T_{N,1,\ell}^{(1)} + \sum_{t=1}^{p-1}(T_{N,1,\ell}^{(t+1)} - T_{N,1,\ell}^{(t)}) + T_{N,2,\ell}^{(p)} + T_{N,3\ell} + T_{N,4\ell}^{(p)}, \qquad (10)$$

where

$$\begin{aligned}
T_{N,1,\ell}^{(t)} = &\left[\sum_{n=1}^{N}\sum_{j=d+1}^{\infty}(K_{[j-1]}(\tilde{\mathbf{X}}_{n,j-1}) - K_{[j]}(\tilde{\mathbf{X}}_{n,j}))\right. \\
& - \sum_{r=1}^{t-1}\sum_{n=1}^{N}\sum_{d+1\leq j_1<\cdots<j_r<\infty}\mathbf{B}_{j_1\cdots j_r} \circ K_{\infty}(\mathbf{0})\prod_{s=1}^{r}\varepsilon_{n+d-j_s} \\
& \left. - \sum_{n=1}^{N}\sum_{d+1\leq j_1<\cdots<j_t<\infty}\mathbf{B}_{j_1\cdots j_t} \circ K_{[j_t]}(\tilde{\mathbf{X}}_{n,j_t})\prod_{s=1}^{t}\varepsilon_{n+d+1-j_s}\right] \\
& - \left[\sum_{n=1}^{N}\sum_{j=d+1}^{\infty}(K_{[j-1]}(\tilde{\mathbf{X}}_{n,j-1,\ell}) - K_{[j]}(\tilde{\mathbf{X}}_{n,j,\ell}))\right.
\end{aligned}$$

$$-\sum_{r=1}^{t-1}\sum_{n=1}^{N}\sum_{d+1\le j_1<\cdots<j_r\le \ell}\mathbf{B}_{j_1\cdots j_r}\circ K_{[\ell]}(\mathbf{0})\prod_{s=1}^{r}\varepsilon_{n+d-j_s}$$

$$-\sum_{n=1}^{N}\sum_{d+1\le j_1<\cdots<j_t\le\ell}b_{j_1\cdots j_t}K_{[j_t]}(\tilde{\mathbf{X}}_{n,j_t,\ell})\prod_{s=1}^{t}\varepsilon_{n+d+1-j_s}\Bigg], \qquad 0\le t\le p,$$

$$T_{N,2,\ell}^{(p)}=\left[\sum_{n=1}^{N}\sum_{d+1\le j_1<\cdots<j_p<\infty}\left(\prod_{s=1}^{t+1}\varepsilon_{n+d-j_s}\right)\right.$$

$$\times\{\mathbf{B}_{j_1\cdots j_p}\circ[K_{[j_p]}(\tilde{\mathbf{X}}_{n,j_p})-K_{[\infty]}(\mathbf{0})]\}\Bigg]$$

$$-\left[\sum_{n=1}^{N}\sum_{d+1\le j_1<\cdots<j_p\le\ell}\left(\prod_{s=1}^{t+1}\varepsilon_{n+d-j_s}\right)\right.$$

$$\times\{\mathbf{B}_{j_1\cdots j_p}\circ[K_{[j_p]}(\tilde{\mathbf{X}}_{n,j_p,\ell})-K_{[\ell]}(\mathbf{0})]\}\Bigg],$$

$$T_{N,3,\ell}=\left[\sum_{n=1}^{N}\sum_{j=1}^{d}(K_{[j-1]}(\tilde{\mathbf{X}}_{n,j-1})-K_{[j]}(\tilde{\mathbf{X}}_{n,j}))\right]$$

$$-\left[\sum_{n=1}^{N}\sum_{j=1}^{d}(K_{[j-1]}(\tilde{\mathbf{X}}_{n,j-1,\ell})-K_{[j]}(\tilde{\mathbf{X}}_{n,j,\ell}))\right],$$

$$T_{N,4,\ell}^{(p)}=-\sum_{r=1}^{p}\sum_{n=1}^{N}\sum_{\substack{1=j_1<\cdots<j_r\\ j_r\ge\ell+1}}\mathbf{B}_{j_1\cdots j_r}\circ K_{\infty}(\mathbf{0})\prod_{s=1}^{r}\varepsilon_{n+d-j_s}$$

$$+\sum_{r=1}^{p}\sum_{n=1}^{N}\sum_{1=j_1<\cdots<j_r\le\ell}[\mathbf{B}_{j_1\cdots j_r}\circ(K_{[\ell]}(\mathbf{0})-K_{[\infty]}(\mathbf{0}))]\left(\prod_{s=1}^{r}\varepsilon_{n+d-j_s}\right).$$

By applying the same martingale decomposition technique used for (8) to $\mathbf{B}_{j_1\cdots j_t}\circ(K_{[j_t]}(\tilde{\mathbf{X}}_{n,j_t})-K_{[\infty]}(\mathbf{0}))$, we have

$$\mathbf{B}_{j_1\cdots j_t}\circ(K_{[j_t]}(\tilde{\mathbf{X}}_{n,j_t})-K_{[\infty]}(\mathbf{0}))$$
$$=\mathbf{B}_{j_1\cdots j_t}\circ\left\{\sum_{i=1}^{\infty}[K_{[j_t+i-1]}(\tilde{\mathbf{X}}_{n,j_t+i-1})-K_{[j_t+i]}(\tilde{\mathbf{X}}_{n,j_t+i})]\right\},$$

which implies that

$$\sum_{d+1\leq j_1<\cdots<j_t<\infty}\left(\prod_{s=1}^{t}\varepsilon_{n+d-j_s}\right)\mathbf{B}_{j_1\cdots j_t}\circ(K_{[j_t]}(\tilde{\mathbf{X}}_{n,j_t})-K_{[\infty]}(\mathbf{0})) \quad (11)$$
$$=\sum_{d+1\leq j_1<\cdots<j_{t+1}<\infty}\left(\prod_{s=1}^{t}\varepsilon_{n+d-j_s}\right)\mathbf{B}_{j_1\cdots j_t}\circ(K_{[j_{t+1}-1]}(\tilde{\mathbf{X}}_{n,j_{t+1}-1})-K_{[j_{t+1}]}(\tilde{\mathbf{X}}_{n,j_{t+1}})).$$

Similarly, repeating the same procedure, but replacing $\infty$ by $\ell$, we have

$$\sum_{d+1\leq j_1<\cdots<j_t\leq \ell}\left(\prod_{s=1}^{t}\varepsilon_{n+d-j_s}\right)\mathbf{B}_{j_1\cdots j_t}\circ(K_{[j_t]}(\tilde{\mathbf{X}}_{n,j_t,\ell})-K_{[\ell]}(\mathbf{0})) \quad (12)$$
$$=\sum_{d+1\leq j_1<\cdots<j_{t+1}\leq \ell}\left(\prod_{s=1}^{t}\varepsilon_{n+d-j_s}\right)\mathbf{B}_{j_1\cdots j_t}\circ(K_{[j_{t+1}-1]}(\tilde{\mathbf{X}}_{n,j_{t+1}-1,\ell})-K_{[j_{t+1}]}(\tilde{\mathbf{X}}_{n,j_{t+1},\ell})).$$

With the help of (11) and (12), we can express $T_{N,1,\ell}^{(t+1)} - T_{N,1,\ell}^{(t)}$ as

$$T_{N,1,\ell}^{(t+1)}-T_{N,1,\ell}^{(t)}=\sum_{n=1}^{N}\sum_{d+1\leq j_1<\cdots<j_{t+1}<\infty}\left(\prod_{s=1}^{t}\varepsilon_{n+d-j_s}\right)\mathbf{B}_{j_1\cdots j_t}\circ L_{n,j_{t+1},\ell},$$

where

$$L_{n,j,\ell}=\left[K_{[j-1]}(\tilde{\mathbf{X}}_{n,j-1})-K_{[j]}(\tilde{\mathbf{X}}_{n,j})-\varepsilon_{n+d-j}\left(\sum_{u=1}^{d+1}\mathbf{A}_{j,u}\frac{\partial}{\partial x_u}\right)\circ K_{[j]}(\tilde{\mathbf{X}}_{n,j})\right]$$
$$-\left[K_{[j-1]}(\tilde{\mathbf{X}}_{n,j-1,\ell})-K_{[j]}(\tilde{\mathbf{X}}_{n,j,\ell})-\varepsilon_{n+d-j}\left(\sum_{u=1}^{d+1}\mathbf{A}_{j,u}\frac{\partial}{\partial x_u}\right)\circ K_{[j]}(\tilde{\mathbf{X}}_{n,j,\ell})\right]$$
$$\times I(j\leq \ell).$$

We also write

$$T_{N,1,\ell}=\sum_{n=1}^{N}\sum_{j=d+1}^{\infty}L_{n,j,\ell},$$

$$T_{N,2,\ell}=\sum_{n=1}^{N}\sum_{d+1\leq j_1<\cdots<j_p<\infty}\left(\prod_{s=1}^{p}\varepsilon_{n+d-j_s}\right)\mathbf{B}_{j_1\cdots j_p}\circ M_{n,j_p,\ell},$$

$$T_{N,3,\ell}=\sum_{n=1}^{N}\sum_{j=1}^{d}P_{n,j,\ell},$$

where
$$M_{n,j,\ell} = [K_{[j]}(\tilde{\mathbf{X}}_{n,j}) - K_{[\infty]}(\mathbf{0})] - [K_{[j]}(\tilde{\mathbf{X}}_{n,j,\ell}) - K_{[\ell]}(\mathbf{0})]I(j \leq \ell)$$

and
$$P_{n,j,\ell} = [K_{[j-1]}(\tilde{\mathbf{X}}_{n,j-1}) - K_{[j]}(\tilde{\mathbf{X}}_{n,j})] - [K_{[j-1]}(\tilde{\mathbf{X}}_{n,j-1,\ell}) - K_{[j]}(\tilde{\mathbf{X}}_{n,j,\ell})]I(j \leq \ell).$$

By repeating the same argument used in (9), we have

$$\operatorname{cov}(L_{n,j,\ell}, L_{n',j',\ell}) = 0 \qquad \text{if } n - j \neq n' - j', \tag{13}$$

$$\operatorname{cov}\left(\left(\prod_{s=1}^{t} \varepsilon_{n+d-j_s}\right) L_{n,j_{t+1},\ell}, \left(\prod_{s=1}^{t} \varepsilon_{n'+d-j'_s}\right) L_{n',j'_{t+1},\ell}\right) = 0 \tag{14}$$
$$\text{if } n - j_s \neq n' - j'_s \text{ for some } 1 \leq s \leq t+1,$$

$$\operatorname{cov}\left(\left(\prod_{s=1}^{p} \varepsilon_{n+d-j_s}\right) M_{n,j_p,\ell}, \left(\prod_{s=1}^{p} \varepsilon_{n'+d-j'_s}\right) M_{n',j'_p,\ell}\right) = 0 \tag{15}$$
$$\text{if } n - j_s \neq n' - j'_s \text{ for some } 1 \leq s \leq p,$$

$$\operatorname{cov}(P_{n,j,\ell}, P_{n',j',\ell}) = 0 \qquad \text{if } n - j \neq n' - j' \tag{16}$$

and, for any pair of real numbers $C_1$ and $C_2$,

$$\operatorname{cov}\left(C_1 \prod_{s=1}^{r_1} \varepsilon_{n-j_s}, C_2 \prod_{s=1}^{r_2} \varepsilon_{n'-j'_s}\right) = 0 \qquad \text{if } r_1 \neq r_2, \text{ or if } r_1 = r_2 = r, \tag{17}$$
$$\text{but } n - j_s \neq n' - j'_s \text{ for some } 1 \leq s \leq r.$$

In order to estimate the growth rate of $\operatorname{var}(Q_{N,p} - Q_{N,p,\ell})$, we also need to compute the non-zero covariances for $L_{n,j,\ell}$, $M_{n,j,\ell}$ and $P_{n,j,\ell}$. In the following, the results of Lemma 4.1 of the next section are used to bound those covariances. Setting $n - j = n' - j'$, we obtain

$$|\operatorname{cov}(L_{n,j,\ell}, L_{n',j',\ell})| \leq C\left(\sum_{m \geq \ell+1} a_m^2\right) a_j^2 a_{j'}^2 \tag{18}$$

by the second result of Lemma 4.1(i) and

$$|\operatorname{cov}(M_{n,j_p,\ell}, M_{n',j'_p,\ell})| \leq C \sum_{m \geq \ell} a_m^2, \qquad \ell \geq 1, \tag{19}$$

by Lemma 4.1(ii). In addition, by (2) and Jensen's inequality,

$$\mathrm{E} P_{n,j,\ell}^2 \leq C \sup_{1 \leq j \leq \ell} \mathrm{E}[K_{[j]}(\tilde{\mathbf{X}}_{1,j}) - K_{[j]}(\tilde{\mathbf{X}}_{1,j,\ell})]^2 \leq C \mathrm{E}[K(\mathbf{X}_1) - K(\mathbf{X}_{1,\ell})]^2. \qquad (20)$$

Combining the orthogonality properties (13)–(17) and the bounds (18)–(20), we can argue the same way as in Ho and Hsing ([14], proofs of Theorems 3.1 and 3.2) and obtain, as $\ell \to \infty$ and uniformly for $N$,

$$N^{-1} \mathrm{var}(T_{N,1,\ell}^{(1)}) = O(\ell^{-(2\beta-1)}) \qquad (\text{by (18)}),$$

$$N^{-1} \mathrm{var}(T_{N,1,\ell}^{(t+1)} - T_{N,1,\ell}^{(t)}) = O(\ell^{-(2\beta-1)}) \qquad (\text{by (18)}),$$

$$N^{-1} \mathrm{var}(T_{N,2,\ell}^{(p)}) = O(\ell^{-((p+1)(2\beta-1)-1)}) \qquad (\text{by (19)})$$

and

$$N^{-1} \mathrm{var}(T_{N,3,\ell}) = O(\ell^{-\alpha_1}) \qquad (\text{by (20) and (4)}),$$

where the bound $O(\ell^{-\alpha_1})$ should be $o(1)$ if, instead of (4), the weaker condition (C2) is assumed. In addition,

$$N^{-1} \mathrm{var}(T_{N,4,\ell}^{(p)}) = O(\ell^{-(2\beta-1)}),$$

which follows from

$$(\mathbf{B}_{j_1 \cdots j_r} \circ [K_{[\infty]}(\mathbf{0}) - K_{[\ell]}(\mathbf{0})])^2 \leq \mathrm{E}(\mathbf{B}_{j_1 \cdots j_r} \circ [K_{[j]}(\tilde{\mathbf{X}}_{1,j}) - K_{[j]}(\tilde{\mathbf{X}}_{1,j,\ell})])^2$$

$$= O\left(\sum_{i \geq \ell} a_i^2\right).$$

The above inequality, which holds uniformly for all of the configurations of $\{j_1, \ldots, j_r\}$, is implied by equation (1) and the first result of part (i) of Lemma 4.1 Consequently, as $N, \ell \to \infty$, $N^{-1} \mathrm{var}(Q_{N,p} - Q_{N,p,\ell})$ is $o(1)$ under (C2) and

$$N^{-1} \mathrm{var}(Q_{N,p} - Q_{N,p,\ell}) = O(\ell^{-\min\{\alpha_1, 2\beta-1, (p+1)(2\beta-1)-1\}}) \qquad (21)$$

uniformly over $N$, if the stronger condition (4) is assumed.

**Proof of Theorem 1.** We prove (5) by the blocking method (Bernstein [2]). Let $A_N, B_N$, and $\ell_N$ be three increasing sequence of positive integers which satisfy

(1) $A_N + B_N < N$ and $A_N + B_N = o(N)$;
(2) $B_N = o(A_N)$ and $\ell_N = \lceil cB_N \rceil, 0 < c < 1$ ($\lceil \cdot \rceil$ is the greatest integer symbol).

Here, $A_N$ and $B_N$, as in the standard setting of the blocking method, are the size of each block and the length between two adjacent blocks, respectively. The exact values of these three sequences $A_N, B_N$ and $\ell_N$ are to be specified later. Recall $T_n(p, \ell)$ defined in (7). Define $k_N = N/(A_N + B_N)$ and, for $s = 1, 2, \ldots, \lceil k_N \rceil$,

$$m_{N,s}(p,\ell_N) = \sum_{n=(s-1)(A_N+B_N)+1}^{sA_N+(s-1)B_N} T_n(p,\ell_N),$$

$$b_{N,s}(p,\ell_N) = \sum_{n=sA_N+(s-1)B_N+1}^{s(A_N+B_N)} T_n(p,\ell_N),$$

$$M_{N,p} = \sum_{s=1}^{\lceil k_N \rceil} m_{N,s}(p,\ell_N) \quad \text{and} \quad B_{N,p} = \sum_{s=1}^{\lceil k_N \rceil} b_{N,s}(p,\ell_N)$$

and

$$R_{N,p} = Q_{N,p,\ell_N} - (M_{N,p} + B_{N,p}).$$

Set $A_N = \lceil N^a \rceil$ and $B_N = \lceil N^b \rceil$, $1 > a > b > 0$. Then $\ell_N = \lceil cN^b \rceil$, $0 < c < 1$. For large $N$, since $T_n(p,\ell_N), 1 \leq n \leq N$, are $\ell_N$-dependent, $\{m_{N,s}(p,\ell_N), s = 1, \ldots, \lceil k_N \rceil\}$ form an i.i.d. sequence. Then, given the result of Lemma 4.2, we can apply the Berry–Esseen theorem (cf. Theorem 7.4.1 of Chung [5]) to the double array $\{m_{N,s}(p,\ell_N), N \geq 1, s = 1, \ldots, \lceil k_N \rceil\}$ to obtain

$$\sup_x |P(N^{-1/2} M_{N,p} \leq x) - \Phi(\sigma^{-1} x)| = O(N^{-(1-a)/2}). \tag{22}$$

We assume for the moment that $N^{-1/2} Q_{N,p}$ is asymptotically normal with positive variance $\sigma^2$. The role of asymptotic normality will become evident later. Set $\delta_N = N^{-\delta}, \delta > 0$ and $\triangle_{N,p} = N^{-1/2}(Q_{N,p} - M_{N,p})$. Using the Petrov inequality,

$$|P(U+V \leq x) - \Phi(x)|$$
$$\leq \sup_x |P(U \leq x) - \Phi(x)| + \frac{\varepsilon}{\sqrt{2\pi}} + P(|V| \geq \varepsilon)$$

(Petrov [17]), where $U$ and $V$ are any two random variables and $\varepsilon$ is any positive real number, we have

$$\sup_x |P(N^{-1/2} Q_{N,p} \leq x) - \Phi(\sigma^{-1} x)|$$
$$\leq \sup_x |P(N^{-1/2} M_{N,p} \leq x) - \Phi(\sigma^{-1} x)| + \sigma^{-2}(\sqrt{2\pi})^{-1} \delta_N + P(|\triangle_{N,p}| > \delta_N)$$
$$\equiv I_1 + I_2 + I_3,$$

where $I_1 = O(N^{-(1-a)/2})$ as already shown in (22), $I_2 = O(N^{-\delta})$ by definition, and $I_3$ can furthermore be bounded by, using Chebyshev's inequality,

$$I_3 \leq 9 \delta_N^{-2} N^{-1} [\mathrm{E}(Q_{N,p} - Q_{N,p,\ell_N})^2 + \mathrm{E} B_{N,p}^2 + \mathrm{E} R_{N,p}^2]$$
$$\equiv I_{3,1} + I_{3,2} + I_{3,3}.$$

From (21),

$$I_{3,1} = O(N^{-(bQ'-2\delta)}),$$

where $Q' = \alpha_1 \wedge (2\beta - 1) \wedge ((p+1)(2\beta - 1) - 1)$. Because $\{b_{N,s}(p, \ell_N), s = 1, \ldots, \lceil k_N \rceil\}$ also forms an i.i.d. sequence, we can, as justified by Lemma 4.2, apply the Berry–Esseen theorem and obtain

$$I_{3,2} = O(N^{-(a-b-2\delta)}).$$

It is easy to see that

$$I_{3,3} = O(N^{-(1-a-2\delta)}).$$

Combining the above five rates that separately dominate $I_1, I_2, I_{3.1}, I_{3.2}, I_{3.3}$ gives that the optimal choice of $\delta$ will be the largest value that satisfies the following inequalities:

$$1 - a \geq 2\delta, \qquad bQ' - 2\delta \geq \delta, \qquad a - b - 2\delta \geq \delta, \qquad 1 - a \geq 3\delta.$$

After some elementary algebra, we get $\delta = Q'/[3(2Q'+1)]$. It is immediate from $I_3 = o(1)$ that (3) holds. The proof is completed. □

**Proof of Theorem 2.** The proof is similar to that of Theorem 1 with $p = 0$. Recall that

$$P_{n,j,\ell} = [K_{[j-1]}(\tilde{\mathbf{X}}_{n,j-1}) - K_{[j]}(\tilde{\mathbf{X}}_{n,j})]$$
$$- [K_{(j-1)^{(1)},j-1}(\tilde{\mathbf{X}}_{n,j-1,\ell}) - K_{[j]}(\tilde{\mathbf{X}}_{n,j,\ell})]I(j \leq \ell).$$

Then

$$\text{var}(Q_{N,0} - Q_{N,0,\ell}) \leq R_{N,1,\ell} + R_{N,2,\ell} + R_{N,3,\ell},$$

where, with $j' = n' - n + j$,

$$R_{N,1,\ell} = 8 \sum_{n=1}^{N} \sum_{j=d+1}^{\ell} \sum_{n'=n}^{n+\ell-j} \text{cov}(P_{n,j,\ell}, P_{n',j',\ell}),$$

$$R_{N,2,\ell} = 8 \sum_{n=1}^{N} \sum_{n'=n}^{N} \sum_{j=\ell+1}^{\infty} \text{cov}(P_{n,j,\ell}, P_{n',j',\ell}),$$

$$R_{N,3,\ell} = 4 \sum_{n=1}^{N} \sum_{j=1}^{d} \sum_{n'=n}^{n+d-j} \text{cov}(P_{n,\ell}, P_{n',\ell}).$$

By (iii) of Lemma 4.1 and the fact that the $|a_j|$ are summable,

$$N^{-1}\text{var}(R_{N,1,\ell}) \leq C\left(\sum_{i=\ell+1}^{\infty} a_i^2\right)\left(\sum_{j=1}^{\infty} |a_j|\right)^2 = O(\ell^{-(2\beta-1)})$$

and

$$N^{-1}\mathrm{var}(R_{N,2,\ell}) \leq C\left(\sum_{j=\ell+1}^{\infty}|a_j|\right)^2 = O(\ell^{-2(\beta-1)}).$$

Following the same argument as in proving $N^{-1}\mathrm{var}(T_{N,3,\ell}) = O(\ell^{-\alpha_1})$, we have

$$\mathrm{var}(R_{N,3,\ell}) = O(\ell^{-\alpha_1}).$$

Hence,

$$\mathrm{var}(Q_{N,0} - Q_{N,0,\ell}) = O(\ell^{-(\alpha_1 \wedge (2\beta-2))}).$$

The rest of the proof is identical to that of Theorem 1, except that the rate is

$$I_{3,1} = O(N^{-(bQ'-2\delta)}) \qquad \text{with } Q' = \alpha_1 \wedge (2\beta-2).$$

This concludes the proof. $\square$

## 4. Technical lemmas

Below are two technical lemmas, Lemmas 4.1 and 4.2, that were used in the preceding section to prove the two main theorems. Lemma 4.1 is the multivariate version of Lemma 6.2 of Ho and Hsing [14]. The proof is omitted since it is similar to that of Ho and Hsing [14] and the main task is to directly apply the regularity conditions (C1) and (C2).

**Lemma 4.1.** *Assume that conditions (C1) and (C2) hold. Let $0 \leq i_1 + \cdots + i_{d+1} \leq J$. Then, for some universal constant $C$,*

(i) *for $j \geq d+1$ and $\ell \geq j$,*

$$\mathrm{E}\left[\frac{\partial^{i_1+\cdots+i_{d+1}}(K_{[j]}(\tilde{\mathbf{X}}_{1,j}) - K_{[j]}(\tilde{\mathbf{X}}_{1,j,\ell}))}{\partial x_1^{i_1}\cdots\partial x_{d+1}^{i_{d+1}}}\right]^2 \leq C\sum_{m=\ell+1}^{\infty} a_m^2$$

*and*

$$\mathrm{E}\left[\frac{\partial^{i_1+\cdots+i_{d+1}}L_{n,j,\ell}}{\partial x_1^{i_1}\cdots\partial x_{d+1}^{i_{d+1}}}\right]^2 \leq C\left(\sum_{m=j}^{\infty} a_m^2\right)\left(\sum_{m=\ell+1}^{\infty} a_m^2\right);$$

(ii) *for $j \geq d+1$ and $\ell \geq j$,*

$$\mathrm{E}\left[\frac{\partial^{i_1+\cdots+i_{d+1}}M_{n,j,\ell}}{\partial x_1^{i_1}\cdots\partial x_{d+1}^{i_{d+1}}}\right]^2 \leq C\left(\sum_{m=j}^{\infty} a_m^2\right)\left(\sum_{m=\ell+1}^{\infty} a_m^2\right);$$

(iii) *for $j \geq 2$,*

$$E\{[K_{[j-1]}(\tilde{\mathbf{X}}_{n,j-1}) - K_{[j]}(\tilde{\mathbf{X}}_{n,j})] - [K_{[j-1]}(\tilde{\mathbf{X}}_{n,j-1,\ell}) - K_{[j]}(\tilde{\mathbf{X}}_{n,j,\ell})]I(j \leq \ell)\}^2$$
$$\leq C\left(\sum_{m=\ell+1}^{\infty} a_m^2\right)(a_j^2).$$

In the preceding Lemma 4.1, we use the fact that there exists a constant $C$ such that $\max_{j-d-1 \leq i \leq j+d+1} |a_i| \leq C|a_j|$.

In Lemma 4.2 below, a moment inequality of fourth order for $Q_{h(N),p,\ell_N}$ is established so that the blocking method can be applied. To prove the lemma, the representation used in the previous section for $Q_{h(N),p,\ell_N}$ (see the identity (10)) is needed. For the sake of presentation, we recall it below.

$$\begin{aligned}
Q_{N,p,\ell} &= \sum_{n=1}^{N} K(X_{n,\ell}^*) - \sum_{r=0}^{p} Z_{N,r,\ell} \\
&= \sum_{n=1}^{N} T_n(p,\ell) \\
&= S_{N,1,\ell}^{(1)} + \sum_{t=1}^{p-1}(S_{N,1,\ell}^{(t+1)} - S_{N,1,\ell}^{(t)}) + S_{N,2,\ell}^{(p)} + S_{N,3,\ell} + S_{N,4,\ell}^{(p)},
\end{aligned} \tag{23}$$

where the various $S$ quantities are defined as follows:

$$S_{N,1,\ell}^{(t+1)} - S_{N,1,\ell}^{(t)} = \sum_{n=1}^{N} \sum_{d+1 \leq j_1 < \cdots < j_{t+1} < \infty} \left(\prod_{s=1}^{t} \varepsilon_{n+d-j_s}\right) \mathbf{B}_{j_1 \cdots j_t} \circ L_{n,j_{t+1},\ell}^*$$

and

$$S_{N,1,\ell} = \sum_{n=1}^{N} \sum_{j=d+1}^{\infty} L_{n,j,\ell}^*,$$

where

$$L_{n,j,\ell}^* = \left[K_{[j-1]}(\tilde{\mathbf{X}}_{n,j-1,\ell}) - K_{[j]}(\tilde{\mathbf{X}}_{n,j,\ell})\right.$$
$$\left. - \varepsilon_{n+d-j}\left(\sum_{u=1}^{d+1} \mathbf{A}_{j,u}\frac{\partial}{\partial x_u}\right) \circ K_{[j]}(\tilde{\mathbf{X}}_{n,j,\ell})\right]I(j \leq \ell),$$

$$S_{N,2,\ell} = \sum_{n=1}^{N} \sum_{d+1 \leq j_1 < \cdots < j_p < \infty} \left( \prod_{s=1}^{p} \varepsilon_{n+d-j_s} \right) \mathbf{B}_{j_1 \cdots j_p} \circ M^*_{n,j_p,\ell},$$

$$S_{N,3,\ell} = \sum_{n=1}^{N} \sum_{j=1}^{d} P^*_{n,j,\ell},$$

where

$$M^*_{n,j,\ell} = [K_{[j]}(\tilde{\mathbf{X}}_{n,j,\ell}) - K_{[\ell]}(\mathbf{0})] I(j \leq \ell)$$

and

$$P^*_{n,j,\ell} = [K_{[j-1]}(\tilde{\mathbf{X}}_{n,j-1,\ell}) - K_{[j]}(\tilde{\mathbf{X}}_{n,j,\ell})] I(j \leq \ell),$$

$$S^{(p)}_{N,4,\ell} = -\sum_{r=1}^{p} \sum_{n=1}^{N} \sum_{1=j_1 < \cdots < j_r \leq \ell} [\mathbf{B}_{j_1 \cdots j_r} \circ (K_{[\ell]}(\mathbf{0}) - K_{[\infty]}(\mathbf{0}))] \left( \prod_{s=1}^{r} \varepsilon_{n+d-j_s} \right).$$

**Lemma 4.2.** *If conditions (C1) and (C2) hold, then*

$$\sup_{N} \mathrm{E} \left( h(N)^{-1/2} \sum_{n=1}^{h(N)} T_n(p, \ell_N) \right)^{4} < \infty, \tag{24}$$

*where $\ell_N$ and $h(N)$ are increasing sequences of positive integers less than $N$ which diverge to $+\infty$.*

**Proof.** In view of the representation (23), it suffices to show that (24) holds for each of

$$S^{(1)}_{h(N),1,\ell_N}, \qquad (S^{(t+1)}_{h(N),1,\ell_N} - S^{(t)}_{h(N),1,\ell_N}), \qquad t = 1, \ldots, p-1,$$

$$S^{(p)}_{h(N),2,\ell_N}, \qquad S_{h(N),3,\ell_N} \quad \text{and} \quad S^{(p)}_{h(N),4,\ell_N}.$$

We only prove the case of $S^{(p)}_{h(N),2,\ell_N}$ since the other cases are similar or simpler (for the case of $S^{(1)}_{h(N),1,\ell}$). Because the index $s$ in $S^{(p)}_{h(N),2,\ell_N}$ only has a finite range, we can assume that it is fixed. Using the orthogonality property given in (15), we get

$$\mathrm{E} \left\{ \sum_{n=1}^{h(N)} \sum_{d+1 \leq j_1 < \cdots < j_p \leq \ell} \left( \prod_{s=1}^{p} \varepsilon_{n+d-j_s} \right) \mathbf{B}_{j_1 \cdots j_p} \circ M_{n,j_p,\ell} \right\}^{4}$$

$$\leq C \left\{ h(N) \sum_{k_1,k_2,k_3=0}^{h(N)-1} \sum_{d+1 \leq j_1 < \cdots < j_p < \infty} \left[ \prod_{u=1}^{p} (j_u(j_u+k_1)(j_u+k_2)(j_u+k_3))^{-\beta} \right] \right. \tag{25}$$

$$\times (j_p(j_p + k_1)(j_p + k_2)(j_p + k_3))^{-(\beta - 1/2)}$$
$$+ \left[ h(N) \sum_{k=0}^{h(N)-1} \sum_{d+1 \le j_1 < \cdots < j_p < \infty} \left( \prod_{u=1}^{p} (j_u(j_u + k))^{-\beta} \right) (j_p(j_p + k))^{-(\beta - 1/2)} \right]^2 \right\},$$

for some constant $C$ independent of $N$ and $\ell$. To derive the right-side of (25), we make use of part (ii) of Lemma 4.1, as well as the fact that the coefficient $\mathbf{A}_{j_1, u_1} \cdots \mathbf{A}_{j_p, u_p}$ in $\mathbf{B}_{j_1 \cdots j_r}$ is bounded above by $C(j_1 \cdots j_p)^{-\beta}$. Because of the assumption $(p+1)(2\beta - 1) > 1$,

$$\sum_{k=0}^{h(N)-1} \sum_{2 \le j_1 < \cdots < j_p < \infty} \left( \prod_{u=1}^{p} (j_u(j_u + k))^{-\beta} \right) (j_p(j_p + k))^{-(\beta - 1/2)} = O(1). \qquad (26)$$

Hence, the second part inside the braces on the right-hand side of (25) is of order $O(h^2(N))$. Similarly, the first part of the right-hand side of (25) is bounded by

$$h^{-2}(N) \left\{ h(N) \sum_{k=0}^{h(N)-1} \sum_{2 \le j_1 < \cdots < j_p < \infty} \left[ \prod_{u=1}^{p} (j_u(j_u + k))^{-\beta} \right] [j_p(j_p + k)]^{-(\beta - 1/2)} \right\}^3 \qquad (27)$$
$$= O(h(N)).$$

Combining (26) and (27) gives $\mathrm{E}(S_{h(N),2,\ell_N}^{(p)})^4 = O(h^2(N))$. Hence, (24) follows. The proof is thus completed. $\square$

## Acknowledgements

The authors are grateful for the referee's careful reading of the manuscript and many valuable suggestions.